\documentclass{elsart}
\usepackage{amsfonts}
\usepackage[dvips]{graphicx}

\usepackage{amsmath}

\usepackage{amssymb}

\begin{document}
\def\diam{{\rm diam}}
\def\ess{\rm ess}
\def\p{{\mathbb P}}
\def\ep{\varepsilon}
\def\P{{\rm Po}}
\def\cf{{\mathcal{ F}}}
\def\cl{{\mathcal{ L}}}
\def\e{{\mathbb E}}
\def\v{{\mathbb V}}
\def\l{\lambda}
\def\ll{{\ell_n}}
\def\a{{\alpha_n}}
\def\ph{\varphi(n)\sqrt n}
\def\dist{{\rm dist}}
\def\lr{\left(}
\def\rr{\right)}
\def\cd{\cdot}
\def\ts{\thinspace}
\def\lc{\left\{}
\def\rc{\right\}}
\def\qed{\vbox{\hrule\hbox{\vrule\kern3pt\vbox{\kern6pt}\kern3pt\vrule}\hrule}}
\newcommand{\hyp}{\mathcal{H}}
\newcommand{\hypp}{\mathcal{K}}
\newcommand{\Nu}{$\begin{Large}$\nu$\end{Large}$}
\newcommand{\ignore}[1]{}
\begin{frontmatter}



\title{Threshold and Complexity Results for the Cover Pebbling Game}


\author{Anant P. Godbole}
\address{East Tennessee State University} \ead{godbolea@mail.etsu.edu}
\author{Nathaniel G. Watson}
\address{Washington University in St. Louis}
\ead{ngwatson@artsci.wustl.edu}
\author{Carl R. Yerger}
\address{Georgia Institute of Technology}
\ead{cyerger@math.gatech.edu}

\begin{abstract}
Given a configuration of pebbles on the vertices of a graph, a
\emph{pebbling move} is defined by removing two pebbles from some
vertex and placing one pebble on an adjacent vertex. The cover
pebbling number of a graph, $\gamma(G)$, is the smallest number of
pebbles such that through a sequence of pebbling moves, a pebble
can eventually be placed on every vertex simultaneously, no matter
how the pebbles are initially distributed.  We determine Bose
Einstein and Maxwell Boltzmann cover pebbling thresholds for the
complete graph.  Also, we show that the cover pebbling decision
problem is NP-complete.
\end{abstract}

\end{frontmatter}

\section {Games}  There are several popular games that involve the movement
of ``objects" along a graph-like structure.  These include
Mancala, where beads are moved along a bent path, and Peg
Solitaire, where pegs are moved across a triangular grid.  In each
case, some objects are removed from the game after a move is made.
Both games proceed according to a set of rules that may be seen at
{\tt http://www.centralconnector.com/GAMES/mancala.html} and {\tt
http://www.mazeworks.com/peggy/} respectively. Mancala is a
well-defined game between two players, while the solo player in
Peg Solitaire pits herself against ``nature". This paper concerns
the {\it pebbling} and {\it cover pebbling} games, which take
place between a highly intelligent Player 1 and a rather
non-competitive opponent with a limited strategy.  There is a
version of pebbling, called ``pegging", which is far closer to Peg
Solitaire on graphs than are our games; see \cite{moulton}. See
also  \cite{chunggraham} for a chessboard game related to
pebbling.

The focus of our paper, and indeed of all previous research on the
subject, is on deriving conditions under which Player 1 wins the
game, or wins with probability that is asymptotic to one, or wins
with probability that approaches zero as the size of the problem
grows to infinity.    Let us start with some baseline definitions
and previously derived facts.
\section{Preliminaries}
Given a connected graph $G$, distribute $t$ pebbles on its
vertices in some configuration. Specifically, a configuration of
weight $t$ on a graph $G$ is a function $C$ from the vertex set
$V(G)$ to $\mathbb{N} \cup \{0\} $ such that $\sum_{v\in
V(G)}C(v)=t$.  Clearly $C$ represents an arrangement of pebbles on
$V(G)$. If the pebbles are indistinguishable, there are
${{n+t-1}\choose{t}}={{n+t-1}\choose{n-1}}$ configurations of $t$
pebbles on $n$ vertices.  Using quantum mechanical terminology as
in \cite{feller}, we shall call this situation {\it Bose Einstein}
pebbling and posit that the underlying probability distribution is
uniform, i.e. that each of the ${{n+t-1}\choose{n-1}}$
distributions are equally likely -- should the pebbles be thrown
randomly onto the vertices. This is the model studied in
\cite{thresh}.  Now there is no reason to assume, {\it a priori},
that the pebbles are indistinguishable.  Accordingly, if the
pebbles are distinct, we shall refer to our process as {\it
Maxwell Boltzmann} pebbling, in which a random distribution of
pebbles leads to $n^t$ equiprobable configurations.  Maxwell
Boltzmann pebbling does not appear to have been studied more than
peripherally in the literature.

A {\it pebbling move} is defined as the removal of two pebbles
from some vertex and the placement of one of these on an adjacent
vertex. Given an initial configuration, a vertex $v$ is called
{\it reachable} if it is possible to place a pebble on it in
finitely many pebbling moves.  The graph $G$ is said to be {\it
pebbleable} (this is not standard nomenclature) if any of its
vertices can be thus reached. Define the pebbling number $\pi(G)$
to be the minimum number of pebbles that are sufficient to pebble
the graph regardless of the initial configuration. The pebbling
game may thus be described as follows:  Player 2 specifies a
distribution $C$ and a target vertex $v$.  Player 1 wins the game
iff she is able to reach vertex $v$ using a sequence of pebbling
moves.  The pebbling number of $G$ is the smallest number $t_0$ of
pebbles so that Player 1 wins no matter what strategy Player 2
employs.

The origin of pebbling is rather interesting and somewhat
unexpected.  Lagarias and Saks \cite{glennsurvey} were considering
a way to produce an alternative proof to a conjecture of Erd\H{o}s
and Lemke, which Lemke and Kleitman proved in 1989
\cite{kleitman}. It is known that for any set $N = \{n_1, n_2,
\ldots, n_q\}$ of $q$ natural numbers, there is a nonempty index
set $I \subset \{1, \ldots , q\}$ such that $q | \sum_{i \in I}
n_i$. The Erd\H os-Lemke conjecture states that the additional
conclusion $\sum_{i \in I} n_i \leq $lcm $(q,n_1,n_2,\ldots, n_q)$
could also be reached (\cite{glennsurvey}). Unfortunately, Lemke
and Kleitman's argument was detailed and contained a considerable
amount of case analysis. This provoked Lagarias and Saks to invent
graph pebbling as a way to produce a cleaner proof, since such a
proof would follow easily if the pebbling number of the cartesian
product of paths was determined.  This was accomplished in a
landmark paper by Chung \cite{chung}.  One generalization of
pebbling, called $p$-pebbling, was utilized in Chung's proof, and
is defined as the removal of $p$ pebbles from some vertex and the
placement of one pebble on an adjacent vertex. It turns out that a
greedy-like condition, called the numerical pebbling operation for
pebbling paths, can be used to prove Chung's theorem
\cite{glennsurvey}. In fact, one of the lemmas of Chung's proof
actually uses the fact that for any set $N = \{n_1, n_2, \ldots,
n_q\}$ of $q$ natural numbers, there is a nonempty index set $I
\subset \{1, \ldots , q\}$ such that $q | \sum_{i \in I} n_i$.

\medskip

\noindent SPECIAL CASES:  The pebbling number $\pi(P_n)$ of the
path is $2^{n-1}$. Chung \cite{chung} proved that $\pi(Q^d)=2^d$
and $\pi(P_n^m)=2^{(n-1)m}$, where $Q^d$ is the $d$-dimensional
binary cube.    An easy pigeonhole principle argument yields
$\pi(K_n)=n$.  The pebbling number of trees has been determined
see \cite{glennsurvey}.

One of the key conjectures in pebbling, now proved in several
special cases, is due to Graham; its resolution would clearly
generalize Chung's result for $m$-dimensional grids:

\bigskip

\noindent GRAHAM'S CONJECTURE.  {\it The pebbling number of the
cartesian product of two graphs is no more that the product of the
pebbling numbers of the two graphs, i.e.
$$\pi(G\Box H)\le \pi(G)\pi(H).$$}

Structural characteristics of graphs have also been employed to
determine the pebbling number of specific classes of graphs.  For
instance, a graph is said to be {\it Class 0} if $\pi(G) = \vert
G\vert$.  Cubes are of Class 0, as are complete graphs, but what
other families fall in this important class of graphs for which
$\pi$ is as low as it can possibly be?  Here are two answers:  For
graphs of diameter $2$, if $G$ is $3$-connected, that is, the
removal of 2 or fewer vertices does not disconnect the graph, then
$G$ is Class 0 \cite{glennhoch}.  In fact, if we consider
$G(n,p)$, the class of random graphs on $n$ vertices where the
probability of each particular edge being present is a fixed
constant $p\in(0,1)$, then almost all such graphs are Class 0
\cite{glennhoch}. Generalizations of this result to the case where
$p=p_n\to0$ as $n\to\infty$ are also available; see
\cite{glenntrotter}.  Other authors, e.g. \cite{chan}, have
obtained general pebbling bounds, while Bukh \cite{bukh} has
proved almost-tight asymptotic bounds on the pebbling number of
diameter three graphs.

Another aspect of pebbling that has been explored is the random
structure one obtains when placing pebbles randomly on graphs.
Specifically, we seek the probability that a graph $G$ is
pebbleable when $t$ pebbles are placed randomly on it according to
the Bose Einstein or Maxwell Boltzmann scheme.  Numerous {\it
threshold results} have been determined in \cite{thresh} for Bose
Einstein pebbling of families of graphs such as $K_n$, the
complete graph on $n$ vertices; $C_n$, the cycle on $n$ vertices;
stars; wheels; etc.  A threshold result is a theorem of the
following kind:
$$t=t_n\gg a_n\Rightarrow \p(G=G_n\ {\rm is\ pebbleable})\to1\quad(n\to\infty)$$
$$t=t_n\ll b_n\Rightarrow \p(G=G_n\ {\rm is\ pebbleable})\to0\quad(n\to\infty),$$
where we write, for non-negative sequences $c_n$ and $d_n$,
$c_n\gg d_n$ (or $d_n\ll c_n$) if $c_n/d_n\to\infty$ as
$n\to\infty$.  Of course, we have reason to feel particularly
gratified if we can show that $a_n=b_n$ in a result of this genre.
For the families of complete graphs, wheels and stars, for
example, we know \cite{thresh} that $a_n=b_n=\sqrt n$.  In many
cases, however, the analysis is quite delicate; see \cite{godbole}
and \cite{path} for some of the issues involved in finding the
pebbling threshold for a family as basic as $P_n$, the path on $n$
vertices.  The fundamental reference \cite{brightwell} contains
general results on the existence of sharp pebbling thresholds for
families of graphs.

A detailed survey of graph pebbling has been presented by Hurlbert
\cite{glennsurvey}, and it would probably not be an
oversimplification to state that most results available to date
fall in four broad categories:  finding pebbling numbers for
classes of graphs; addressing the issue of when a family of graphs
is of class 0; pinpointing graph pebbling thresholds; and seeking
to understand the complexity issues in graph pebbling
\cite{glenncomplex}. A survey of open problems in graph pebbling
may be found on Glenn Hurlbert's website; see {\tt
http://math.la.asu.edu/ $\sim$hurlbert/HurlPebb.ppt}.

The above mini-survey on pebbling notwithstanding, we focus in
this paper on a {\it variant} of pebbling called cover pebbling,
first discussed by Crull et al \cite{crull}.  For reasons that
will become obvious, we focus only on analogs of the last two of
the four general directions mentioned above.

The {\it cover pebbling number} $\lambda(G)$ is defined as the
minimum number of pebbles required such that it is possible, given
any initial configuration of at least $\lambda(G)$ pebbles on $G$,
to make a series of pebbling moves that {\it simultaneously}
reaches {\it each} vertex of $G$.  A configuration is said to be
{\it cover solvable} if it is possible to place a pebble on every
vertex of $G$ starting with that configuration. Various results on
cover pebbling have been determined. For instance, we now know
\cite{crull} that $\l(K_n)=2n-1; \l(P_n)=2^n-1$; and that for
trees $T_n$,
\begin{equation}
\l(T_n)=\max_{v\in V(T_n)}\sum_{u\in
V(T_n)}2^{\dist(u,v)}.\end{equation} Likewise, it was shown in
\cite{munyon} that $\l(Q^d)=3^d$ and in \cite{nate} that
$\l(K_{r_1,\ldots,r_m})=4r_1+2r_2+\ldots+2r_m-3$, where $r_1\ge
r_2\ge \ldots \ge r_m$.  The above examples reveal that for these
special classes of graphs at any rate, the cover pebbling number
equals the ``stacking number", or, put another way, the worst
possible distribution of pebbles consists of placing all the
pebbles on a single vertex.  The intuition built by computing the
value of the cover pebbling number for the families $K_n$, $P_n$,
and $T_n$ in \cite{crull} led to Open Question No.~10 in
\cite{crull}, which was christened the {\it Stacking Conjecture}
by students at the Summer 2004 East Tennessee State University
REU.  In an exciting summer development, participants Annalies
Vuong and Ian Wyckoff \cite{ian} were able to prove the

\bigskip

\noindent STACKING THEOREM:  {\it For any connected graph $G$,
$$\l(G)=\max_{v\in V(G)}\sum_{u\in V(G)}2^{\dist(u,v)},$$}

\bigskip
\noindent thereby proving that (1) holds for all graphs.  (The
stacking theorem was independently proved soon after by
Sjostr\"and \cite{sjostrand}.) In fact, the key result in Vuong
and Wyckoff's paper \cite{ian} is really a sufficient condition
for a distribution to be cover solvable, so further investigations
in the theory of (cover) pebbling might soon veer, we speculate,
in a fifth general direction, namely a study of which
distributions are (cover) solvable and which are not.  Indeed,
such a research thrust would be most consistent with our
description of pebbling as a {\it game}, and is addressed in
Section 5 of this paper.

\section{Maxwell Boltzmann Cover Pebbling Threshold for $K_n$}  It is evident that $n$ is the smallest number of pebbles that might suffice to cover pebble $K_n$ -- in the unlikely event that they happen to be distributed one apiece on the vertices.  On the other hand, we know that $2n-1$ pebbles always suffice, since $\l(K_n)=2n-1$.  We seek a sharp cover pebbling threshold that is somewhere in between these two extremes, when distinguishable pebbles are thrown onto the $n$ vertices of $K_n$ according to the Maxwell Boltzmann scheme.  To explain {\it why} there is a dramatic increase in the probability of the pebbleability of $K_n$ at $t=(1.5238.....)n$, we first prove an important auxiliary result that gives necessary and sufficient conditions for a configuration to be cover solvable.  Such results are not easy to come by, as we will further see in Section 5.
\subsection{Necessary and Sufficient Conditions for Cover Solvability of $K_n$}
Let $X=X_{n,t}$ be the number of vertices on which an odd number
of pebbles are placed. We will often refer to $X$ as the number of
{\it odd stacks}.
\begin{lem}
A configuration of $t$ pebbles on the $n$ vertices of $K_n$ is
cover solvable if and only if
\begin{equation} X+t\ge 2n.
\end{equation}
\end{lem}

\medskip

\noindent{\bf Proof 1} Let $C$ be any cover solvable
configuration.  This implies that after a sequence of pebbling
moves, each previously uncovered vertex has two pebbles associated
with it -- one on it, and one that was removed from the game.
Likewise, any vertex that previously had a non-zero and  even
number of pebbles on it must have at least two pebbles still on
it, and any  vertex that previously had an odd number of pebbles
on it must now have at least one pebble on it.  Thus we must have
$t\ge2E+X$ (where $E$ represents the number of vertices on which
an even number (including 0) of pebbles were initially placed),
which is a condition that is easily seen to be equivalent to (2).

Conversely, if $C$ is not cover solvable, then we must, after a
series of pebbling moves, reach a point where there are still some
uncovered vertices; where each initially covered vertex in the set
$E$ now has exactly two pebbles on it; and each vertex in the set
$X$ has exactly one pebble on it.  It follows that
$t<2E+X$.\hfill\qed

\medskip

\noindent{\bf Proof 2}  Given a configuration $C$, let $Y_i,0\le
i\le t$, be the number of vertices with $i$ pebbles. Now a vertex
can cover exactly $i$ others if and only if it has either $2i+1$
or $2i+2$ pebbles on it.  It follows that $C$ is cover solvable
iff
\[\sum_{i\ge0}i(Y_{2i+1}+Y_{2i+2})\ge Y_0,\]
or, iff
\[\sum_{i\ge0}(2i+1)Y_{2i+1}+\sum_{i\ge0}(2i+2)Y_{2i+2}\ge2Y_0+\sum_{i\ge0}Y_{2i+1}+2\sum_{i\ge0}Y_{2i+2},\]
i.e., iff $t\ge 2E+X$.\hfill\qed

Armed with Lemma 1, we now provide the heuristic reason why we
believe there is a sharp cover pebbling threshold at
$t=(1.5238....)n$.  Given a random variable $X$ with expected
value $\e(X)$, we will say that $X$ is {\it sharply concentrated}
around $\e(X)$ if $X\sim \e(X)$ with high probability (w.h.p.)
Assuming therefore that $X\sim\e(X)$ w.h.p., it makes sense to
speculate that $K_n$ is cover solvable with high probability
whenever $\e(X)\ge 2n-t$.  But $X=\sum_{j=1}^nI_j,$ where $I_j=1$
(resp.~0) if there is an odd (resp.~even) stack on vertex $j$, so
that linearity of expectation yields
\begin{eqnarray}\e(X)&=&n\p(I_1=1)\nonumber\\
&=&n\sum_{j\ {\rm odd}}{t \choose j} \lr {1\over n} \rr^j\lr 1-{1\over n}\rr^{t-j}\nonumber\\
&=& {n\over2}\lr 1- \lr1-{2\over n}\rr^t\rr.
\end{eqnarray}
Thus $\e(X)\ge 2n-t$ iff
\begin{equation}t-{n\over 2}\lr1-{2\over n}\rr^t\ge{{3n}\over{2}},\end{equation}
and, parametrizing by setting $t=An$, we see that (4) holds iff
\begin{equation}
A-{1\over 2}\lr1-{2\over n}\rr^{An}\ge{3\over 2}.
\end{equation}
 Since $(1-2/n)^n\sim e^{-2}$ we see from (5) that a reasonable guess for a threshold value of $t$ is $A_0n$ where $A_0$ is the solution of
\[A-{1\over 2}\exp\{-2A\}={3\over 2},\]
or $A_0$=1.5238.......  We next make this hunch more precise.
\subsection{Main Result}  Various tools are used to establish concentration of measure results.  Some of the more sophisticated techniques employed are the martingale method, a.k.a.~Azuma's inequality (proved independently and a few years earlier by W.~Hoeffding), and Talagrand's isoperimetric inequalities in product spaces; see \cite{steele} for an exposition of both.  In this section, however, we establish our main result by using a baseline technique, {\it viz.} Tchebychev's inequality, also known in probabilistic combinatorics circles as the ``second moment method."  Theorem 2 below thus yields the ``correct" result, but with a sub-optimal rate of convergence due to the crudeness of the method employed.  Tchebychev's inequality states that for any random variable $X$ and for any $\l>0$,
\[\p\lr\vert X-\e(X)\vert\ge\l\rr\le{{\v(X)}\over{\l^2}},\]
where $\v(Z)$ denotes the variance of $Z$.  We start by computing
the variance of $X$; $\e(X)$ has already been evaluated in (3).
We have
\begin{eqnarray}
\v(X)&=&\v\lr\sum_{j=1}^nI_j\rr
\nonumber\\&=&\sum_{j=1}^n\v(I_j)+\sum_{i\ne j}\lr\e(I_iI_j)-\e(I_i)\e(I_j)\rr\\
&=&n\p(I_1=1)(1-\p(I_1=1))+n(n-1)\lc\p(I_1I_2=1)-\p^2(I_1=1)\rc,\nonumber
\end{eqnarray}
so we focus on computing $\p(I_1I_2=1)$, i.e., the probability
that both vertex 1 and vertex 2 have odd stacks of pebbles.  It is
not too hard to verify that
\begin{eqnarray*}
\p(I_1I_2=1)&=&\sum_{r,s\ {\rm odd}}{{t}\choose{r,s,t-r-s}}\lr{1\over n}\rr^r\lr{1\over n}\rr^s\lr1-{2\over n}\rr^{t-r-s}\nonumber\\
&=&{1\over4}\lr\Sigma_1-\Sigma_2-\Sigma_3+\Sigma_4\rr,
\end{eqnarray*}
where
\[\Sigma_1=\sum_{r,s}{{t}\choose{r,s,t-r-s}}\lr{1\over n}\rr^r\lr{1\over n}\rr^s\lr1-{2\over n}\rr^{t-r-s},\]
\[\Sigma_2=\sum_{r,s}{{t}\choose{r,s,t-r-s}}\lr{1\over n}\rr^r\lr-{1\over n}\rr^s\lr1-{2\over n}\rr^{t-r-s},\]
\[\Sigma_3=\sum_{r,s}{{t}\choose{r,s,t-r-s}}\lr-{1\over n}\rr^r\lr{1\over n}\rr^s\lr1-{2\over n}\rr^{t-r-s},\]
and
\[\Sigma_4=\sum_{r,s}{{t}\choose{r,s,t-r-s}}\lr-{1\over n}\rr^r\lr-{1\over n}\rr^s\lr1-{2\over n}\rr^{t-r-s},\]
so that
\[
\p(I_1I_2=1)={1\over4}\lr1+\lr1-{4\over n}\rr^t-2\lr1-{2\over
n}\rr^t\rr,
\]
and
\begin{eqnarray*}
{\rm Cov}(I_1,I_2)&=&\e(I_1I_2)-\e(I_1)\e(I_2)\\
&=&{1\over 4}\lc\lr1-{4\over n}\rr^t-\lr1-{2\over n}\rr^{2t}\rc.
\end{eqnarray*}
It follows from (6) that
\begin{equation}
\v(X)={n\over 4}\lr1-\lr1-{2\over n}\rr^{2t}\rr+{{n(n-1)}\over
{4}}\lc\lr1-{4\over n}\rr^t-\lr1-{2\over n}\rr^{2t}\rc.
\end{equation}
We are now ready to state
\begin{thm}
Consider $t$ distinct pebbles that are thrown onto the vertices of
the complete graph $K_n$ on $n$ vertices according to the Maxwell
Boltzmann distribution.  Set $A_0=1.5238\ldots$.  Then
\[t= A_0n+\varphi(n){\sqrt{n}}\Rightarrow\p(K_n\ {\rm is\ cover\ solvable})\to1\quad(n\to\infty)\]
and
\[t= A_0n-\varphi(n){\sqrt{n}}\Rightarrow\p(K_n\ {\rm is\ cover\ solvable})\to0\quad(n\to\infty),\]
where $\varphi(n)\to\infty$ is arbitrary.\end{thm}

\medskip

\noindent{\bf Proof}  Assume first that $t=
A_0n+\varphi(n){\sqrt{n}}$.  Then \begin{eqnarray}
&&{}\p(X\ge 2n-t)\nonumber\\
&&=\p\lr X-\e(X)\ge2n-t-{n\over 2}\lr1-\lr1-{2\over n}\rr^t\rr\rr\nonumber\\
&&=\p\lr X-\e(X)\ge{3\over 2}n-A_0n+{n\over 2}\lr1-{2\over n}\rr^{A_0n+\varphi(n)\sqrt{n}}-\varphi(n)\sqrt{n}\rr\nonumber\\
&&=\p\lr X-\e(X)\ge-{n\over2}e^{-2A_0}+{n\over 2}\lr1-{2\over n}\rr^{A_0n+\varphi(n)\sqrt{n}}-\varphi(n)\sqrt{n}\rr\nonumber\\
&&\ge\p\lr X-\e(X)\ge-{n\over2}e^{-2A_0}+{n\over 2}\exp\lc-2A_0-{{2\varphi(n)}\over{\sqrt{n}}}\rc-\varphi(n)\sqrt{n}\rr\nonumber\\
&&=\p\lr X-\e(X)\ge{n\over2}e^{-2A_0}\lc\exp\lc{-{2\varphi(n)}\over{\sqrt{n}}}\rc-1\rc-\varphi(n)\sqrt{n}\rr\nonumber\\
&&\ge\p\lr X-\e(X)\ge{n\over2}e^{-2A_0}\cdot-{{2\varphi(n)}\over{\sqrt{n}}}(1+o(1))-\varphi(n)\sqrt{n}\rr\nonumber\\
&&=\p\lr X-\e(X)\ge-\varphi(n)\sqrt{n}\lr1+e^{-2A_0}\rr(1+o(1))\rr\nonumber\\
&&\ge\p\lr \vert
X-\e(X)\vert\le\varphi(n)\sqrt{n}\lr1+e^{-2A_0}\rr(1+o(1))\rr.
\end{eqnarray}
In (8), the first and second inequalities follow from the facts
that $1-x\le e^{-x}$ and $e^{-x}-1\le-x/(1+x)$ respectively.  We
next analyze the variance as given by (7).  The first component
${n\over 4}(1-(1-{2\over n})^{2t})$ of $\v(X)$ can easily be
verified to be of the form ${n\over4}(1+K(1+o(1)))$ (for some
constant $K$) when $t$ is as specified; the second component
\small{
\[{{n(n-1)}\over {4}}\lc\lr1-{4\over n}\rr^t-\lr1-{2\over n}\rr^{2t}\rc={{n(n-1)}\over {4}}\lc\lr1-{4\over n}\rr^t-\lr1-{4\over n}+{4\over{n^2}}\rr^{t}\rc,\]
}on the other hand, may be bounded using the inequalities
\[t(b-a)a^{t-1}\le b^t-a^t\le t(b-a)b^{t-1}\]
for $b>a$ to yield, for $t$ as above,
\begin{eqnarray*}&&{} {{n(n-1)}\over {4}}\lc\lr1-{4\over
n}+{4\over{n^2}}\rr^{t}-\lr1-{4\over n}\rr^t\rc
\\&& ={{n(n-1)}\over{4}}{{4t}\over{n^2}}\lr1-{4\over
n}+\Theta\lr{1\over
{n^2}}\rr\rr^{t-1}\\&&=\Theta(t)\\&&=\Theta(n).\end{eqnarray*} It
follows that $\v(X)=\Theta(n)$ and thus we have by (8) and
Tchebychev's inequality, with $K$ representing a generic constant,
\[\p(X\ge 2n-t)\ge\p\lr \vert X-\e(X)\vert\le K\cdot\sqrt{n}\varphi(n)\rr\ge1-{{1}\over{K^2\varphi^2(n)}}\to1\quad(n\to\infty),\]
as asserted.

The proof of the second half of the theorem is similar; we bound
above instead of below to get, with $t=
A_0n-\varphi(n){\sqrt{n}}$,
\begin{eqnarray*}
&&{}\p(X\ge 2n-t)\nonumber\\
&&=\p\lr X-\e(X)\ge-{n\over2}e^{-2A_0}+{n\over 2}\lr1-{2\over n}\rr^{A_0n-\varphi(n)\sqrt{n}}+\varphi(n)\sqrt{n}\rr\nonumber\\
&&\le\p\bigg( X-\e(X)\ge-{n\over2}e^{-2A_0} \\ &&\phantom{asd} +{n\over 2}\exp\lc\lr-2A_0+{{2\varphi(n)}\over{\sqrt{n}}}\rr\lr1+o(1)\rr\rc +\varphi(n)\sqrt{n}\bigg)\nonumber \\
&&=\p\lr X-\e(X)\ge{n\over2}e^{-2A_0}\lc\exp\lc{{2\varphi(n)}\over{\sqrt{n}}}(1+o(1))\rc-1\rc  +\varphi(n)\sqrt{n}\rr \nonumber \\
&&\le\p\lr X-\e(X)\ge{n\over2}e^{-2A_0}\cdot{{2\varphi(n)}\over{\sqrt{n}}}(1+o(1))+\varphi(n)\sqrt{n}\rr\nonumber\\
&&\le\p\lr \vert X-\e(X)\vert\ge\varphi(n)\sqrt{n}\lr1+e^{-2A_0}(1+o(1))\rr\rr\\
&&\le\frac{K}{\varphi^2(n)}\to0.
\end{eqnarray*}
This completes the proof. \hfill\qed

\noindent Remarks: Note that $\v(X)=\Theta(n)$ implies that $X$ is
concentrated w.h.p.~in an interval of length  $\Omega(\sqrt{n})$
around $\e(X)=\Theta(n)$.  This is because of the inequality
$\v(X)\le {\rm Range}^2(X)/4$.

\section{Bose Einstein Cover Pebbling Threshold}
\subsection{Exact Distributions}  In the Maxwell Boltzmann scheme, it is extremely difficult (though not impossible) to calculate $\p(X=x)$ exactly.  Moreover, the formula for $\p(X=x)$ obtained thus is quite intractable.  Surprisingly, this is not the case when we consider Bose Einstein statistics; we can derive such a formula by counting the number of configurations of size $t$ on $K_n$ with $x$ odd stacks. Call this value $ \phi (x, t, n)$.  Clearly if $x$ and $t$ have different parity, or $x>t$ or $x>n,$ $\phi(x,t,n)=0.$ Suppose this is not the case. Then, we can find the configurations with $x$ odd stacks of pebbles by placing one pebble on the $x$ vertices that are to have odd stacks on them, and then distributing the remaining $t-x$ pebbles on the $n$ vertices of $G$ in $\frac{t-x}{2}$ indistinguishable pairs. Thus, since the vertices with odd stacks may be chosen in ${{n}\choose{x}}$ ways, we have proved
\begin{prop} If $t$ and $x$ have the same parity, and if $x\le\min\{t,n\}$, then
\[\p(X=x)={{{n \choose x} { \frac{t-x}{2}+n-1 \choose n-1}}\over{{n+t-1}\choose{n-1}}}.\]
\end{prop}
In fact, the above exact distribution was used fruitfully in
\cite{nate2} to prove a weaker version of the main result of the
next section. We choose, however, to prove a cover pebbling
threshold for Bose Einstein pebbling by using more contemporary
probabilistic tools.
\subsection{Polya Sampling and Azuma's Inequality Yield Dividends}
There is a natural and sequential probabilistic process associated
with Maxwell Boltzmann pebbling.  We simply take $t$ pebbles
(balls) and throw them one by one onto (into) $n$ vertices (urns)
in the ``natural" way that inspires many elementary problems in
discrete probability texts.  By contrast, the ``global"
Bose-Einsteinian positioning of $t$ indistinct balls into $n$
distinct urns -- so that we obtain ${n+t-1 \choose n-1}$
equiprobable configurations -- does not {\it appear} to have a
sequential process associated with it.  {\it But it does.}  We
first rephrase the problem -- not as one associated with throwing
balls {\it into} boxes but, conversely, as a sampling problem,
i.e., drawing balls {\it from} boxes.  In this light, Maxwell
Boltzmann pebbling consists of drawing $t$ balls ``with
replacement" from a box containing one ball of each of $n$ colors,
with the understanding that the number of balls of color $j$ drawn
in the altered model equals the number of pebbles that are tossed
onto vertex $j$ {\it \`a la} the balls-in-boxes model.  {\it Bose
Einstein pebbling can be recast in a similar fashion, but one
needs to employ a process called P\'olya sampling.}  P\'olya
sampling (or the P\'olya urn model) is described in \cite{feller}
as a means of modelling contagious diseases and takes place as
follows:  Initially the urn contains one ball of each of $n$
colors.  After each draw, the selected ball is replaced together
with another ball of the same color.  In this mode of sampling, we
lose the independence inherent to the with-replacement procedure,
and, as a matter of fact, the selection process is not even
Markovian -- but are able to ``see" the sequentiality that will be
critical in the sequel.  As before, the number of times that color
$j$ is drawn can be set to equal the number of pebbles on vertex
$j$, but do these two procedures yield the same probability model?
We claim so, and here is a proof of this rather well-known fact:
\begin{lem}
Let $X_j$ be the number of times the color $j$ is drawn among the
$t$ draws.  Then for any $x_1,x_2,\ldots,x_n$ with $\sum x_j=t$,
$$\p(X_1 =x_1, X_2 =x_2,\ldots, X_n =x_n) =
{{1}\over{{{n+t-1}\choose{t}}}}$$
\end{lem}

\medskip

\noindent {\bf Proof} First let us find the probability that the
stated outcome appears in the following order:  color 1 is first
drawn $x_1$ times, then color 2 is drawn $x_2$ times, etc., until
we finally draw color $n$ the last $x_n$ times.  Call this ordered
event $A$.  It is easily seen that
\begin{eqnarray*}\p(A) &=& {{(1\cdot 2\cdot\ts\cdots \ts\cdot x_1)(1\cd 2\cd \ts \cdots\ts \cd
x_2)\cdots(1\cd 2\cd \ts\cdots \ts\cd
x_n)}\over{n(n+1)\cdots(n+t-1)}}\\ &=& {{x_1!x_2!\cdots
x_n!}\over{n(n+1)\cdots(n+t-1)}}.\end{eqnarray*} Moreover, a
little reflection reveals that this probability is the same
regardless of the order in which the balls are drawn. The total
number of ways of ordering our configuration turns out to be
${{t!}/{x_1!x_2!\cdots x_n!}}$.  Thus, we can write the
probability of obtaining a given configuration as
\begin{eqnarray*}\p(X_1 =x_1, X_2 =x_2,\ldots, X_n =x_n)& =& {{t!}\over{x_1!x_2!\cdots
x_n!}}\cdot {{x_1!x_2!\cdots x_n!}\over{n(n+1)\cdots(n+t-1)}}\\&
=& {{1}\over{{{n+t-1}\choose{t}}}}.\end{eqnarray*}This concludes
the proof.  \hfill\qed

We next provide the reader with background concerning the Azuma
martingale inequality.  Consider the following generic set up:
$(\Omega,\cf, \p)$ is a probability space and $(Y_n)$ a sequence
of random variables on $(\Omega,\cf, \p)$ that may not, in
general, be independent. In our case, the $Y_j$s are the sequence
of draws made in the P\'olya urn scheme associated with the random
cover pebbling problem.  Let $X=X_t=X(Y_1,\dots,Y_t)$ be an
objective function (in our case $X$ is the number of odd stacks of
pebbles), and consider the filtration (sequence of sigma algebras)
$(\cf_n)$, $\cf_0=\{\emptyset,\Omega\}$,
$\cf_i=\sigma(Y_1,\dots,Y_i)$. Let $\e_iX$ denote the conditional
expectation of $X$ with respect to $\cf_i$ (with $\e=\e_0$) and
set $d_i =\e_iX-\e_{i-1}X$. Then it is well known that
$(d_i,\cf_i)$ is a martingale difference sequence, and that we
have $X-\e( X)=\sum_{i=1}^td_i$. A key method used towards gaining
an understanding of the concentration of $X$ around $\e(X)$ is the
method of bounded differences, also known as the Azuma
(1967)-Hoeffding (1963) inequality \cite{steele}:
\begin{lem} (Azuma-Hoeffding) For all $\l>0$, \begin{equation}\p(\vert X-\e(
X)\vert\ge\l)\le2\exp\lc-{\l^2\over2\sum\|d_i\|_{\infty}^2}\rc.\end{equation}\end{lem}
\noindent where $\|Z\|_\infty={(\ess)\sup}|Z(\omega)|$.  (9) may
be made further applicable as follows:  Letting $(Y_i^*)$ be an
independent copy of $(Y_i)$, we have
$$\e_{i-1}X(Y_1,\ldots,Y_t)=\e_iX(Y_1,\ldots,Y_{i-1},Y_i^*,
Y_{i+1}\ldots,Y_t),$$ so that $d_i$ can be written as a single
conditional expectation as follows:
$$d_i=\e_i\big(
X(Y_1,\ldots,Y_t)-X(Y_1,\ldots,Y_{i-1},Y_i^*,
Y_{i+1}\ldots,Y_t)\big),$$ and thus
\begin{eqnarray}\|d_i\|_{\infty}&=&\|\e_i\{
X(Y_1,\ldots,Y_t)-X(Y_1,\ldots,Y_{i-1},Y_i^*,
Y_{i+1}\ldots,Y_t)\}\|_{\infty}\nonumber\\
&\le& \|X(Y_1,\ldots,Y_t)-X(Y_1,\ldots,Y_{i-1},Y_i^*,
Y_{i+1}\ldots,Y_t)\|_{\infty}.
\end{eqnarray}
The philosophy behind the method of bounded differences is thus
that the change in the value of $X$ resulting from a change in a
single input variable should be small.  Moreover the bound in (10)
shows us, after a moment's reflection, that we have, for our
problem, $\|d_i\|_\infty\le2$. Furthermore, Lemma 5 and (10) yield
the following concentration for the number $X$ of odd stacks in
the Bose-Einstein scheme:
$$\p\lr\vert X-\e(X)\vert\ge\l\rr\le2\exp\lc-{{\l^2}\over{8t}}\rc,\eqno({\rm Az})$$
so that $X$ is concentrated in an interval of length
$\sqrt{n}\varphi(n)$ around $\e(X)$ whenever $t\sim Kn$.  We are
now ready to prove the main result of this section -- one that
features the golden ratio $\gamma=(1+\sqrt{5})/2$:
\begin{thm}
Consider $t$ indistinguishable pebbles that are placed on the
vertices of the complete graph $K_n$ according to the Bose
Einstein distribution. Then, with $\gamma$ representing the golden
ratio $(1+\sqrt{5})/2$,
\[t= \gamma n+\varphi(n){\sqrt{n}}\Rightarrow\p(K_n\ {\rm is\ cover\ solvable})\to1\quad(n\to\infty)\]
and
\[t= \gamma n-\varphi(n){\sqrt{n}}\Rightarrow\p(K_n\ {\rm is\ cover\ solvable})\to0\quad(n\to\infty),\]
where $\varphi(n)\to\infty$ is arbitrary.\end{thm}

\noindent{\bf Proof}  We start by establishing tight bounds on
$\e(X)$; unfortunately, the sum
\[\e(X)=n\p(I_1=1)\]
does not appear to be tractable.  Consider first an upper bound on
$\p(I_1=1)$:  We have, with $\a$ and $\ll$ representing,
respectively, generic functions of the form $O(1/n)$ and $O(\log
n/n)$ whose exact form may vary from line to line,
\begin{eqnarray}
\p(I_1=1)&=&\p({\rm vertex\ 1\ has\ an\ odd\ number\ of\ pebbles})\nonumber\\
&=&\sum_{j=1}^{``\infty"}\frac{{{n+t-2j-1}\choose{n-2}}}{{{n+t-1}\choose{n-1}}}\nonumber\\
&=&\sum_{j=1}^{\infty}\frac{n-1}{n+t-1}\frac{t}{n+t-2}\frac{t-1}{n+t-3}\frac{t-2}{n+t-4}\ldots\frac{t-2j+2}{n+t-2j}\nonumber\\
&\le&\frac{n}{n+t}\frac{t}{n+t}(1+\a)\sum_{j=1}^\infty \frac{t-1}{n+t-3}\frac{t-2}{n+t-4}\ldots\frac{t-2j+2}{n+t-2j}\nonumber\\
&\le&\frac{nt}{(n+t)^2}(1+\a)\sum_{j=1}^\infty\lr\frac{t-1}{n+t-3}\rr^{2j-2}\nonumber\\
&=&(1+\a)\frac{nt}{(n+t)^2}\frac{1}{1-\lr\frac{t-1}{n+t-3}\rr^2}\nonumber\\
&=&(1+\a)\frac{nt}{(n+t)^2}\frac{1}{1-\lr\frac{t}{n+t}\rr^2}\nonumber\\
&=&(1+\a)\frac{t}{n+2t}.
\end{eqnarray}
Next note that for $R$ to be determined,
\begin{eqnarray}
\p(I_1=1)&\ge&\sum_{j=1}^R \frac{n-1}{n+t-1}\frac{t}{n+t-2}\frac{t-1}{n+t-3}\frac{t-2}{n+t-4}\ldots\frac{t-2j+2}{n+t-2j}\nonumber\\
&\ge& \frac{nt}{(n+t)^2}(1+\a)\sum_{j=1}^R\lr\frac{t-2j+2}{n+t-2j}\rr^{2j-2}\nonumber\\
&\ge& \frac{nt}{(n+t)^2}(1+\a)\sum_{j=1}^R\lr\frac{t-2R+2}{n+t-2R}\rr^{2j-2}\nonumber\\
&=&\frac{nt}{(n+t)^2}(1+\a)\frac{1-\lr\frac{t-2R+2}{n+t-2R}\rr^{2R}}{1-\lr\frac{t-2R+2}{n+t-2R}\rr^2}.
\end{eqnarray}
We now pick $R$ in (12) so that
\[1-\lr\frac{t-2R+2}{n+t-2R}\rr^{2R}=1-\a;\]
this may be done, e.g., with $R=K\log n$, since in the range of
$t$'s that we are dealing with (namely $t=\gamma
n\pm\varphi(n)\sqrt n$), we have
\[1-\lr\frac{t-2R+2}{n+t-2R}\rr^{2R}=1-\lr\frac{\gamma\lr1+O\lr\frac{\varphi(n)}{\sqrt n}\rr\rr}{(1+\gamma)\lr1+O\lr\frac{\varphi(n)}{\sqrt n}\rr\rr}\rr^{2R}=1+\a\]
if $R$ is an appropriate multiple of $\log n$.  Equation (12) thus
yields
\begin{eqnarray}
\p(I_1=1)&\ge&\frac{nt}{(n+t)^2}(1+\a)\frac{(n+t-2R)^2}{(n-2)(n+2t-4R+2)}\nonumber\\
&=&(1+\a)(1+\ll)\frac{t}{n+2t}\nonumber\\
&=&(1+\ll)\frac{t}{n+2t}.
\end{eqnarray}
Equations (11) and (13) thus give
\begin{equation}
(1+\ll)\frac{nt}{n+2t}\le\e(X)\le(1+\a)\frac{nt}{n+2t},
\end{equation}
and we are ready to derive our cover pebbling threshold.  We first
set $t=\gamma n+\varphi(n)\sqrt n$ to get, for some (and this is
important) {\it positive} $A$,
\begin{eqnarray}&&\p(X\ge 2n-t)\nonumber\\
&&\ge \p\lr X-\e(X)\ge 2n-t-\frac{nt}{n+2t}(1+\ll)\rr\nonumber\\
&&=\p\lr X-\e(X)\ge(2-\gamma)n-\varphi(n)\sqrt n-\frac{n(\gamma n+\ph)}{n+2(\gamma n+\ph)}(1+\ll)\rr\nonumber\\
&&=\p\lr X-\e(X)\ge(2-\gamma) n-\ph-\frac{\gamma n^2\lr1+\frac{\varphi(n)}{\gamma\sqrt n}\rr}{n(1+2\gamma)\lr1+\frac{2\varphi(n)}{(1+2\gamma)\sqrt n}\rr}(1+\ll)\rr\nonumber\\
&&\ge\p\bigg( X-\e(X)\ge(2-\gamma)n-\ph  \nonumber\\ && \phantom{A}-\frac{\gamma n}{(1+2\gamma)}\lr1 +A\lr\frac{\varphi(n)}{\sqrt n}\rr\rr(1+\ll)\bigg).\nonumber\\
\end{eqnarray}
Now $(2-\gamma)n-{\gamma n}/{(1+2\gamma)}=0$, so (15) yields
\begin{eqnarray}
\p(X\ge 2n-t)&\ge&\p\lr X-\e(X)\ge -\ph-\frac{A\gamma}{1+2\gamma}\ph+O(\log n)\rr\nonumber\\
&\ge&\p(X-\e(X)\ge-B\ph+O(\log n))\nonumber\\
&\ge&\p(\vert X-\e(X)\vert\le B\ph+O(\log n))\nonumber\\
&\ge&1-2\exp\lc-\frac{(B\ph)+O(\log n))^2}{8t}\rc\nonumber\\
&\to&0,
\end{eqnarray}
by Azuma's inequality as given by Equation (Az).

Conversely, for $t=\gamma n-\ph$, we have for some $C>0$,
\begin{eqnarray}
&&\p(X\ge 2n-t)\nonumber\\
&&\le\p\lr X-\e(X)\ge(2-\gamma)n+\ph-\frac{n(\gamma n-\ph)}{n+2(\gamma n-\ph)}(1+\a)\rr\nonumber\\
&&=\p\lr X-\e(X)\ge(2-\gamma)n+\ph-\frac{\gamma n^2\lr1-\frac{\varphi(n)}{\gamma\sqrt n}\rr(1+\a)}{(1+2\gamma)n\lr1-\frac{2\varphi(n)}{(1+2\gamma)\sqrt n}\rr}\rr\nonumber\\
&&\le\p\lr X-\e(X)\ge(2-\gamma)n+\ph-\frac{n\gamma}{1+2\gamma}\lr1-\frac{C\varphi(n)}{\sqrt n}\rr(1+\a)\rr\nonumber\\
&&=\p(X-\e(X)\ge D\ph+O(1))\nonumber\\
&&\le\p(\vert X-\e(X)\vert\ge D\ph+O(1))\nonumber\\
&&\to0,
\end{eqnarray}
again by (Az).  Equations (16) and (17) complete the
proof.\hfill\qed

\section{NP-Completeness of the Cover Pebbling Problem} One of the
obvious open problems that can be formulated as a result of our
work is the following:  What are cover pebbling thresholds for
families of graphs other than $K_n$?  It would certainly advance
the theory of cover pebbling if one could uncover a host of
results similar, e.g. to Theorems 2 and 6. Such results would
provide a nice complement to those in \cite{thresh}.  Our results
in this section show, however, that this task might not be as easy
as one might imagine. Necessary and sufficient conditions for the
cover solvability of a graph are likely to be complicated, and the
best hope might thus be to establish necessary conditions and
sufficient conditions that are not too far apart.

The normal way to formalize the concept of the difficulty of a
problem is to use the concept of computational complexity.
Formally, we imagine a {\it decision  problem} to be a set of
infinite strings of characters (like data represented by bits in a
computer.) A decision problem is said to {\it accept} a string if
this set contains the string. Usually, we look for the best
possible asymptotic upper bound (in terms of the length of the
string) for the number of steps the fastest possible algorithm
takes to determine whether a given string is in the set.
Informally, we think of decision problems being yes-no questions
about a property of some class of finite mathematical structures
(graphs, integer matrices, etc.) and we ask how fast it is
possible to correctly determine the yes or no answer in terms of
the size of the input.

For instance, some problems can be solved by an algorithm which
takes only a number of steps which is bounded by a polynomial in
the size of the input, while others take at least an exponential
amount of time to solve. The former class of decision problems is
called $P$ for ``polynomial." The class $N\!P,$ for
``nondeterministic polynomial" is a bit more complicated; roughly
speaking it is the set of decision problems for which a ``yes"
answer can be ``checked'' in polynomial time, given an appropriate
piece of information. That is, if we call the class of inputs to
the decision problem $X,$ and the class of inputs which the
decision problem accepts $X',$ there exists a class $Y$ of objects
(called the {\it certificates}) and a function $A \ : \ X \times Y
\rightarrow \{0,1\}$ which is computable in polynomial time, such
that for any instance $x \in X$ of the decision problem, there
exists a $y \in Y$ such that $A(x,y)=1$ if and only if $x \in X'.$
For instance, the decision problem which asks whether a given
number is composite is easily seen to be in $N\!P,$ because the
composite numbers are exactly those with nontrivial divisor, and,
given two numbers, it is easy to determine by division whether one
is a divisor of the other. Also, any problem in $P$ is also in
$N\!P,$ because any polynomial-time method of solving a problem is
trivially also a polynomial-time method of verifying a yes answer.
However, it is a celebrated open problem if the converse also
holds and $P=N\!P.$

Within $N\!P,$ there is a class of problems, called the {\it
$N\!P$-complete} problems, which is thought of being a set of
problems which are at least as hard as any other problem in
$N\!P.$ This is because any instance $x$ of a problem $D$ in
$N\!P$ can be translated by a polynomial-time algorithm to an
instance $x'$ of any $N\!P$-complete problem $D'$ such that $x'$
is accepted by $D'$ if and only if $x$ is accepted by $D.$
Therefore, if we could solve any $N\!P$-complete problem in
polynomial time, we could solve any problem in $N\!P$ in
polynomial time by translating it to an instance of this problem.
Thus, the question of whether $P=N\!P$ reduces to the question of
whether any particular $N\!P$-complete problem can be solved in
polynomial time.

We now show that the problem which asks if a configuration of
pebbles on a graph is cover solvable is $N\!P$-complete. It is
worth noting that most complexity theorists speculate that $P
\not= N\!P,$ and therefore, when a problem is classified as
$N\!P$-complete, it is usually thought of as evidence of its
difficulty.  See \cite{Garey} for a comprehensive theory of
$N\!P$-completeness, and Watson \cite{nate3} and Milans et
al.~\cite{milans} for a more general exposition on the complexity
of cover pebbling.

\begin{thm}
Let $G$ be a graph, and $C$ a configuration on $G.$  Let $|G|=m$
and label the vertices of $G$ as $v_1, v_2,\ldots v_m$. Then $C$
is cover solvable if and only if there exist integers $n_{ij}\geq
0$  with $1\leq i,j \leq m$ and $n_{ij}=0$ and $n_{ji}=0$ whenever
$\{v_i,v_j\} \notin E(G) $ such that for all $1\leq k \leq m,$
$$C(v_k)+\sum_{l=1}^m n_{lk}-2\sum_{l=1}^m n_{kl} \geq 1.$$
\end{thm}

\medskip

\noindent{\bf Proof} First, suppose $C$ is cover solvable. Then
find some sequence of pebbling moves which cover solves $C.$ Let
$n_{ij}$ be the total number of pebbling moves from $v_i$ to $v_j$
in this sequence. Then after all the moves, there are exactly
$$C(v_k)+\sum_{l=1}^m n_{lk}-2\sum_{l=1}^m n_{kl}$$ pebbles left
on $v_k$, which is always at least $1$ because of the fact that
this sequence of moves cover solves $C.$

Conversely, suppose such numbers $n_{ij}$ exist.  This means that
there {\it does} exist a sequence of moves that solves $C$, with
$n_{ij}$ moves being made from $v_i$ to $v_j$, but possibly with
some illegal ``negative pebbling" along the way.  We show,
however, that for each $i,j$ it is possible to {\it legally} make
$n_{ij}$ moves from $v_i$ to $v_j$; since for each $k$,
$$C(v_k)+\sum_{l=1}^m n_{lk}-2\sum_{l=1}^m n_{kl} \geq 1,$$ which
leads to a cover solution of $C$.  The main question thus is:  In
what order do we make these moves?  We proceed in any arbitrary
fashion, continuing to make pebbling moves as long as there exist
vertices $v_{i'}$ and $v_{j'}$ such that less than $n_{i'j'}$
moves from  $v_{i'}$ to $v_{j'}$ have already been made and there
are at least two pebbles on $v_{i'}.$ If no such pair $\{v_{i'},
v_{j'}\}$ exists, then for each $(i,j)$, either $n_{ij}$ moves
have been made from $i$ to $j$ or else there is at most 1 pebble
on vertex $v_i$.  Let $C'$ be the configuration left on $G$ after
these moves and $S$ be the set of $v_i \in G$  for which the total
number of moves from $v_i$ is less than  $\sum_{l=1}^m n_{il}.$

If $S=\emptyset$ then clearly for every $1\leq i,j\leq m$ we have
made $n_{ij}$ moves from $v_i$ to $v_j$ and thus, for every $k$
there are  $C(v_k)+\sum_{l=1}^m n_{lk}-2\sum_{l=1}^m n_{kl}\geq 1$
pebbles on $v_k,$ so $C'(v) \geq 1$ for all $v \in G$ and we have
cover solved $C.$

If  $S \not= \emptyset$ then consider the total number of moves
that remain to be made from a vertex of $S$ (the total $\sum_{i\in
S} \sum_{l=1}^m n_{il}$ minus the number of moves that have
already been made from vertices in $S.$) By the definition of $S,$
this total is at least $|S|,$ since at least one move remains to
be made from every vertex in $S.$ The total is exactly $|S|$ only
if exactly one move remains to be made from each vertex. Also,
$C'(v)\leq 1$ for all $v \in S,$ for a total of at most $|S|$
pebbles. Consider the remaining moves, each of which must
originate from $S.$ Each of these moves, if executed, would remove
one pebble from $S$ if they also end at a vertex in $S,$ and two
if they end at a vertex outside of $S.$ Thus each move must both
begin and end in $S,$ for otherwise $S$ would be left with a
negative total number of pebbles at the {\it end} of the pebbling
sequence, which is impossible. Even if all moves begin and end in
$S$, however, we end up with at most $0$ pebbles on $S$ -- which
too is impossible since we have assumed that the moves
cover-solve, and so there must be one pebble on each vertex of S.
Thus we must have $S=\emptyset$ and we have cover solved $C$.
\hfill\qed

\begin{cor}
The cover solvability decision problem which accepts pairs $\{ G,
C\}$ if and only if $G$ is a graph and $C$ is a configuration
which is cover solvable on $G$ is in $N\!P.$
\end{cor}

\medskip

\noindent{\bf Proof} The above theorem gives the appropriate
certificate of cover solvability, any list of integers $n_{ij}$
which satisfy the equation in Theorem 7. Indeed, Theorem 7 shows
that cover pebbling is equivalent to a special case of the
$N\!P$-complete problem of integer programming, which asks, given
an $n \times m$ integer matrix $A$ and an $n$-dimensional integer
vector $b$ if there exists an $m$-dimensional integer vector $x$
such that $Ax \geq b,$ holds componentwise. Having reduced cover
solvability to a special case of this $N\!P$ problem, we know that
cover solvability is also in $N\!P$.

We now pause to point out another corollary which will be needed
later and which is interesting in its own right:

\begin{cor} Let $G$ be a graph, $C$ a configuration on $G.$ If the sequence of pebbling moves $Q=(q_1, q_2 \ldots q_k) $ cover solves $C,$ and it is possible to make the sequence $Q'=(q_{i_1}, q_{i_2} , \ldots q_{i_l})  $ of moves (with $1 \leq i_j \leq k$ for all $j$ but with no particular requirement on the order of the $i_j,$) then the configuration $C'$ obtained from $C$ after the moves $ (q_{i_1}, q_{i_2} , \ldots q_{i_l})  $ is cover solvable.
\end{cor}

\medskip

\noindent{\bf Proof} The point of this corollary is that the order
of our pebbling moves can't matter. To show this, we simply note
that if it were possible to somehow execute the remaining moves
from $Q$ which are not in $Q',$ they would solve $C'.$ By Theorem
7, it is thus possible to solve $C'.$ \hfill\qed

Now we turn our attention to showing that the cover solvability
decision problem is $N\!P$-hard, that is, that any instance of any
problem in $N\!P$ can be translated to an instance of cover
solvability in polynomial time. The usual method of showing that a
problem $A$ is $N\!P$-hard is to find an $N\!P$-complete problem
$B$ for which any instance of $B$ can be translated into an
instance of $A$ in polynomial time. Then for any instance of any
problem in $N\!P$ we can translate it in polynomial time to an
instance of $B,$ then translate this instance into an instance of
$A.$ For cover solvability, we will use a known $N\!P$-complete
problem known as ``exact cover by 4-sets." Indeed, the
corresponding and seemingly simpler problem of perfect cover by
3-sets is also $N\!P$-complete, but for our purposes, the 4-set
problem is more useful.

\begin{thm} {\it (Karp) \cite{Karp}}
Let the {\it exact cover by 4-sets problem} be the decision
problem which takes as input a set $S$ with $4n$ elements and a
class $A$ of at least $n$ 4-element subsets of $S,$ accepting such
a pair if there exists an $A' \subseteq A$ such that $A'$ is a
class of disjoint subsets which make a partition of $S,$ that is
they are $n$ subsets containing every element of $S.$ This problem
is $N\!P$-complete.
\end{thm}

\begin{thm}
The cover solvability decision problem is $N\!P$-complete
\end{thm}

\medskip

\begin{figure}[htb]
\unitlength 1mm
\begin{center}
\begin{picture}(80,100)
\put(40,0){\circle*{3}} \put(40,10){\circle*{3}}
\put(40,25){\circle*{3}} \put(30,25){\circle*{3}}
\put(50,25){\circle*{3}}
\put(40,10){\line(-2,3){10}} 
\put(40,10){\line(0,1){15}} \put(40,10){\line(2,3){10}} 
\put(40,40){\circle*{3}} 
\put(60,40){\circle*{3}}\put(20,40){\circle*{3}}
\put(40,40){\line(0,1){15}} 
\put(60,40){\line(2,3){10}}
\put(20,40){\line(-2,3){10}} 
\put(40,25){\line(0,1){15}} 
\put(50,25){\line(2,3){10}}
\put(30,25){\line(-2,3){10}} 
\put(40,55){\circle*{3}} 
\put(70,55){\circle*{3}} \put(10,55){\circle*{3}}
\multiput(5,70)(10,0){8}{\circle*{3}} \put(4, 74){$t_1$}
\put(14,74){$t_2$}\put(24, 74){$t_3$} \put(34,74){$t_4$} \put(44,
74){$t_5$} \put(54,74){$t_6$}\put(64, 74){$t_7$}
\put(74,74){$t_8$}
\put(40,55){\line(-1,1){15}} \put(40,55){\line(-1,3){5}}
\put(40,55){\line(1,3){5}} \put(40,55){\line(1,1){15}}
\put(10,55){\line(-1,3){5}} \put(10,55){\line(1,3){5}}
\put(10,55){\line(1,1){15}} \put(10,55){\line(5,3){24}}
\put(70,55){\line(1,3){5}} \put(70,55){\line(-1,3){5}}
\put(70,55){\line(-1,1){15}} \put(70,55){\line(-5,3){24}}
\put(14,53){9} \put(44,53){9} \put(74,53){9} \put(24,38){1}
\put(44, 38){1} \put(64, 38){1} \put(34,23){1} \put(44, 23){1}
\put(54, 23){1} \put(3,53){$b_1$}
\put(33,53){$b_2$}\put(63,53){$b_3$} \put(34, 8){$v$} \put(34,
-2){$w$} \put(44,8){2} \put(44,-2){0} \put(40,0){\line(0,1){10}}
\end{picture}
\end{center}
\caption{\label{dia} A graph that corresponds to the exact cover
by four $4$-sets problem, $a_1 = \{s_1,s_2,s_3,s_4\}$, $a_2 =
\{s_3,s_4,s_5,s_6\}$, $a_3 = \{ s_5,s_6,s_7,s_8 \}$.} \label{ex2}
\end{figure}
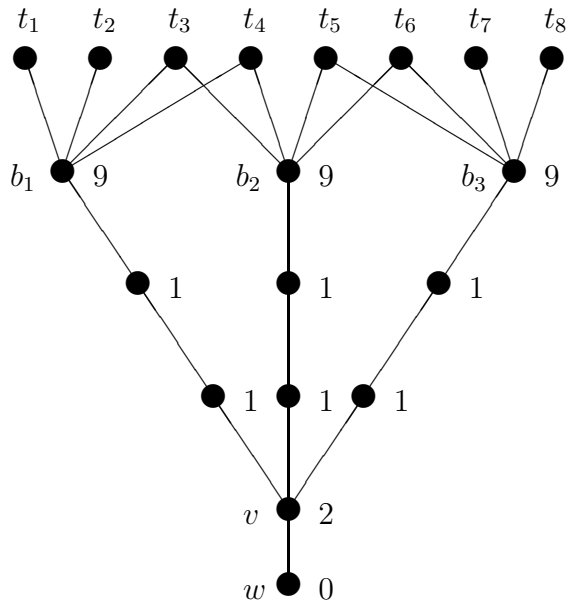

\noindent{\bf Proof} Having shown that this decision problem is in
$N\!P,$ it remains to be shown that it is $N\!P$-hard. Given a set
$S=\{s_1,s_2,\ldots s_{4n}\}$ and a class $A= \{a_1, a_2, \ldots
a_m\}$ of four element subsets of $S,$ that is, an instance of the
exact cover by 4-sets problem,  we construct a graph $G'$ and a
configuration $C'$ on $G'$ in the following manner: We create a
set of vertices $T=\{t_1, t_2, \ldots t_{4n}\}$ which will be
thought of as corresponding to the elements of $S,$ and a set of
vertices $B = \{b_1, b_2, \ldots b_m\}$ which will be thought of
as corresponding to the members of $A.$ Let $C'(t)=0$ for all $t
\in T$ and let $C'(b) =9$ for all $b \in B.$ We make edges between
$B$ and $T$ in the intuitive way, including $\{b_i, t_j\}$ if $t_j
\in b_i.$ Additionally, create a vertex $v$ and a path of length
$m-n$ which has one terminal vertex $v$ and the other called $w.$
Let $C'(v) = 2^{m-n}-(m-n)+1,$ $C'(w)=0$ and $C'(u)=1$ for all $u$
between $v$ and $w$ on the path. Finally, create vertex classes $
B' = \{b_1',b_2' \ldots b_m' \} $ and $B''= \{b_1'',b_2'' \ldots
b_m'' \},$ creating edges $\{b_i, b_i'\},$ $\{b_i',b_i''\}$ and
$\{b_i'', v\}$  for all $i.$ Let $C'(u) = 1 $ for all $u \in B'
\cup B''.$ (Figure \ref{dia}.)

 Clearly, this construction can be made in polynomial time in the size of the pair $\{S, A\}$. Indeed, we have made only $5|A|+1$ vertices and $8|A|-|S|$ edges. In order to finish the proof, we claim that $C'$ is solvable if and only if $A$ contains a perfect cover of $S.$

First suppose that $A$ contains a perfect cover $A' = \{ a_{i_1},
a_{i_2}, \ldots a_{i_n} \} $ of $S.$ Then for each vertex in $B$
which is a $b_{i_j}$ for some $ 1 \leq j \leq n,$ we use 8 of the
pebbles on this $b_{i_j},$ two each to cover the four vertices of
$T$ to which it is adjacent. Because of the fact that $A'$ is a
perfect cover and the way we constructed $G',$ we now have exactly
one pebble on every vertex of $T.$ Furthermore, we have $m-n$
vertices in $B$ which still have 9 pebbles each on them. Because
$v$ is at distance $3$ from each of these vertices, we can use 8
pebbles from each of these vertices to move one pebble each onto
$v$ from these $m-n$ vertices. This leaves $2^{m-n} +1$ pebbles on
$v,$ which is enough to move one pebble onto $w$ while leaving one
pebble on $v.$ After this is done, we have exactly one pebble on
every vertex of $G',$ and we thereby know that $C'$ is solvable.

To show the converse, suppose that $A$ does not contain a perfect
cover of $S.$ Suppose as well that $C'$ is solvable on $G'.$
Clearly, the sequence of pebbling moves which solves $C'$ must
contain (at least) one move to $t$ for every $t \in T.$ Clearly,
each of these moves must originate from $B,$ and no more than 4
can originate from any one vertex of $B.$ Since $A$ does not
contain a perfect cover of $S,$ it cannot be the case that these
moves originate from exactly $n$ vertices in $B.$

We make these $4n$ moves immediately from $C',$ using Corollary 9
to see that the resulting configuration must be solvable (we use
the fact that no more than 4 of these moves can originate from any
one vertex in $B$ to see that it is indeed possible to make these
moves). In addition to the one pebble left on every vertex of $B$
to ensure they remain covered, there are now $8(m-n)$ pebbles on
$B,$ but they are not in $m-n$ groups of $8$ pebbles because the
moves we made originated from more than $n$ vertices of $B.$ In
order to reach $w,$ we clearly need to move $m-n$ pebbles onto $v$
while leaving the rest of the graph covered. Clearly, this is only
possible if all $8(m-n)$ extra pebbles are moved by a path of
length $3$ onto $v.$ However, only one such path is available for
any group of these pebbles.  But any group of less than $8$
pebbles cannot increase the number of pebbles on $v$ by moving
along this path, while leaving the vertices of the path covered.
Since there are not indeed $m-n$ groups of 8 pebbles on $B$ we see
that it is impossible to gather the pebbles necessary to reach $w$
and thus our configuration is not solvable, which is a
contradiction.

\section{Acknowledgements}
This work was done under the supervision of Anant Godbole at the
East Tennessee State University Research Experience for
Undergraduates (REU) site during the summer of 2004, while Watson
and Yerger were undergraduate students at Washingon University in
St.~Louis, and Harvey Mudd College respectively.    Support from
NSF Grant DMS-0139291 is gratefully acknowledged by all three
authors, as is the support offered by fellow REU participants --
``pebblers" and ``non-pebblers" alike.



\end{document}